\documentclass[12pt]{article}%
\usepackage{graphicx}
\usepackage{amsmath}
\usepackage{amsfonts}
\usepackage{amssymb}%
\newtheorem{theorem}{Theorem}

\newtheorem{conjecture}[theorem]{Conjecture}
\newtheorem{corollary}[theorem]{Corollary}

\newtheorem{definition}[theorem]{Definition}

\newtheorem{lemma}[theorem]{Lemma}

\newtheorem{proposition}[theorem]{Proposition}

\begin{document}

\title{On the Conformal Geometry of Transverse Riemann-Lorentz Manifolds}
\author{E. Aguirre, V. Fern\'{a}ndez, J. Lafuente.\\\emph{{\small Dept. Geometr\'{\i}a y Topolog\'{\i}a, Fac. CC. Matem\'{a}ticas,
UCM.}}\\\emph{{\small Plaza de las Ciencias 3, Madrid, Spain.}}}
\maketitle

\begin{abstract}
Physical reasons suggested in \cite{Ha-Ha} for the \emph{Quantum Gravity
Problem} lead us to study \emph{type-changing metrics} on a manifold. The most
interesting cases are \emph{Transverse Riemann-Lorentz Manifolds}. Here we
study the conformal geometry of such manifolds.

\end{abstract}

\section{Preliminaries}

Let $M$ be a connected manifold, $\dim M=m\geq2$, and let $g$ be a symmetric
covariant tensor field of order $2$ on $M$. Assume that the set $\Sigma$ of
points where $g$ degenerates is not empty. Consider $p\in\Sigma$ and $\left(
\mathbb{U},x\right)  $ a coordinate system around $p$. We say that $g$ is a
\emph{transverse type-changing metric} on $p$ if $d_{p}\left(  \det\left(
g_{ab}\right)  \right)  \neq0$ (this condition does not depend on the choice
of the coordinates). We call $\left(  M,g\right)  $ \emph{transverse
type-changing pseudoriemannian manifold} if $g$ is transverse type-changing on
every point of $\Sigma$. In this case, $\Sigma$ is a hypersurface of $M$.
Moreover, at every point of $\Sigma$ there exists a one-dimensional
\emph{radical}, that is the subspace $Rad_{p}\left(  M\right)  $ of $T_{p}M$
which is $g$-\emph{ortogonal to the whole }$T_{p}M$ (and it can be transverse
or tangent to the hypersurface $\Sigma$). The \emph{index} of $g$ is constant
on every connected component of $\mathbb{M}=M-\Sigma$, thus $\mathbb{M}$ is a
union of connected pseudoriemannian manifolds. Locally, $\Sigma$ separates two
pseudoriemannian manifolds whose indices differ in one unit (so we call
$\Sigma$ \emph{transverse type-changing hypersurface}, in particular $\Sigma$
is orientable). The most interesting cases are those in which $\Sigma$
separates a riemannian part from a lorentzian one. We call these cases
\emph{transverse Riemann-Lorentz manifolds}.

Let $\tau\in C^{\infty}(M)$ be such that $\tau\mid_{\Sigma}=0$ and $d\tau
\mid_{\Sigma}\neq0$. We say that (locally, around $\Sigma$) $\tau=0$\emph{\ is
an equation for }$\Sigma$. Given $f\in C^{\infty}(M)$, it holds: $\tau
\mid_{\Sigma}=0\Leftrightarrow f=k\tau$, for some $k\in C^{\infty}(M)$. In
what follows we shall use this fact extensively.

On $\mathbb{M}$ we have naturally defined all the objects associated to
pseudoriemannian geometry, derived from the Levi-Civita connection. In
\cite{Kos85}, \cite{Kos87}, \cite{Kos94}, \cite{Kos97} and \cite{Ag-La}, the
extendibility of geodesics, parallel transport and curvatures have been
studied. Our aim in the present paper is to study the conformal geometry of
transverse Riemann-Lorentz manifolds, including criteria for the extendibility
of the \emph{Weyl conformal curvature}.

Let $\left(  M,g\right)  $ be a transverse Riemann-Lorentz manifold. First of
all, note that we do not have any Levi-Civita connection $\nabla$ defined on
the whole $M$. However we have (\cite{Kos85}) a unique torsion-free metric
\emph{dual connection }
\[
\square:\mathfrak{X}\left(  M\right)  \times\mathfrak{X}\left(  M\right)
\rightarrow\mathfrak{X}^{\ast}\left(  M\right)
\]
on $M$ defined by a \emph{Koszul-like formula}. On $\mathbb{M}$ it holds
$\square_{X}Y\left(  Z\right)  =g\left(  \nabla_{X}Y,Z\right)  $, and thus the
concepts derived from Levi-Civita connection $\nabla$ (on $\mathbb{M}$)
coincide with those derived from the dual connection $\square$.

We say that a vectorfield $R\in\mathfrak{X}\left(  M\right)  $ \emph{is
radical} if $R_{p}\in Rad_{p}\left(  M\right)  -\left\{  0\right\}  $ for all
$p\in\Sigma$. Given a radical vectorfield $R\in\mathfrak{X}\left(  M\right)
$, $\left.  \square_{X}Y\left(  R\right)  \right\vert _{\Sigma}\,$only depends
on $\left.  X\right\vert _{\Sigma}$ and $\left.  Y\right\vert _{\Sigma}$, thus
we obtain the following well-defined map
\[
II^{R}:\mathfrak{X}_{\Sigma}\times\mathfrak{X}_{\Sigma}\rightarrow C^{\infty
}\left(  \Sigma\right)  ,\left(  X,Y\right)  \mapsto\square_{X}Y\left(
R\right)
\]

Note that the $II^{R}$-orthogonal complement to $Rad_{p}\left(  M\right)  $ is
$T_{p}\Sigma$ (\cite{Kos97}, 1(a)), thus $X\in\mathfrak{X}_{\Sigma}$ is
tangent to $\Sigma$ if and only if $II^{R}\left(  X,R\right)  =0$.

Because of the properties of $\square$, the restriction of $II^{R}$ to
vectorfields in $\mathfrak{X}\left(  \Sigma\right)  $ is a well-defined
$\left(  0,2\right)  $ symmetric tensor field $II_{\Sigma}^{R}\in S^{2}\left(
\Sigma\right)  $. Furthermore, since $\square_{X}Y$ is a one-form on $M$ and
the radical is one-dimensional, the condition $II_{\Sigma}^{R}=0$ does not
depend on the radical vectorfield $R$. \emph{A transverse Riemann-Lorentz
manifold is said to be }$II$\emph{-flat} \emph{if }$II_{\Sigma}^{R}=0$\emph{,
for some (and thus, for any) radical vectorfield }$R$. It turns out
(\cite{Kos97} for transverse, \cite{Ag-La} for tangent radical) that $M$ is
$II$-flat if and only if all covariant derivatives $\nabla_{X}Y$, for
$X,Y\in\mathfrak{X}\left(  M\right)  $ tangent to $\Sigma$, smoothly extend to
$M$. Moreover, in that case, $\left.  \nabla_{X}Y\right\vert _{\Sigma}$ only
depends on $\left.  X\right\vert _{\Sigma}$ and $\left.  Y\right\vert
_{\Sigma}$, thus we obtain another well-defined map
\[
III^{R}:\mathfrak{X}\left(  \Sigma\right)  \times\mathfrak{X}\left(
\Sigma\right)  \rightarrow C^{\infty}\left(  \Sigma\right)  ,\left(
X,Y\right)  \mapsto II^{R}\left(  \nabla_{X}Y,R\right)
\]
which is a $\left(  0,2\right)  $ symmetric tensorfield on $\Sigma$. \emph{A
transverse Riemann-Lorentz }$II$\emph{-flat metric is said to be }%
$III$\emph{-flat if }$III^{R}=0$.

If the radical is tangent, $\nabla_{R}R$ becomes transverse (\cite{Ag-La});
therefore, in order that a $II$-flat metric becomes $III$-flat, the radical
must be transverse. And we have the following result (\cite{Kos97}),
concerning the extendibility of curvature tensors:

\begin{theorem}
\label{1.1}The covariant curvature $K$ smoothly extends to $M$ if and only if
the radical is transverse and $g$ is $II$-flat, while the Ricci tensor $Ric$
smoothly extends to $M$ if and only if the radical is transverse and $g$ is
$III$-flat.
\end{theorem}

\section{A Gauss formula for Transverse Riemann-Lorentz Manifolds}

Let $\left(  M,g\right)  $ be a transverse Riemann-Lorentz manifold with
transverse radical.

\begin{lemma}
There exists a unique (canonically defined) radical vectorfield $R$ such that
$II^{R}\left(  R,R\right)  =1$.
\end{lemma}

\textbf{Proof: }Given a radical vectorfield $U$, consider $R=\left(
II^{U}\left(  U,U\right)  \right)  ^{-\frac{1}{3}}\cdot U$, which is a
well-defined radical vectorfield (since the radical is transverse). Thus
$II^{R}\left(  R,R\right)  =1$. Furthermore, if $Z=fR$ is another radical
vectorfield such that $II^{Z}\left(  Z,Z\right)  =1$, then $1=II^{Z}\left(
Z,Z\right)  =f^{3}II^{R}\left(  R,R\right)  =f^{3}$, and consequently $f=1$
$\clubsuit\vspace{0.5cm}$

Suppose that $\left(  M,g\right)  $ is $II$-flat. As we said before, given
$X,Y\in\mathfrak{X}\left(  \Sigma\right)  $, $\nabla_{X}Y$ is well-defined.
Moreover, $\tan\left(  \nabla_{X}Y\right)  :=\nabla_{X}Y-III^{R}\left(
X,Y\right)  \cdot R$ is indeed tangent to $\Sigma$, since%
\[
II^{R}\left(  R,\tan\left(  \nabla_{X}Y\right)  \right)  =III^{R}\left(
X,Y\right)  -III^{R}\left(  X,Y\right)  II^{R}\left(  R,R\right)  =0
\]

\begin{lemma}
\label{2.0.0}If $X,Y\in\mathfrak{X}\left(  \Sigma\right)  $ and $\nabla
^{\Sigma}$ is the Levi-Civita connection of $\left(  \Sigma,g_{\Sigma}\right)
$, it holds
\[
\nabla_{X}Y=\nabla_{X}^{\Sigma}Y+III^{R}\left(  X,Y\right)  \cdot R
\]

\end{lemma}

\textbf{Proof: }Let be $Z\in\mathfrak{X}\left(  \Sigma\right)  $. Since
$\left(  M,g\right)  $ $\,$is $II$-flat, $\nabla_{X}Y$ is well defined and it
must hold $\square_{X}Y\left(  Z\right)  =g\left(  \nabla_{X}Y,Z\right)
=g_{\Sigma}\left(  \tan\left(  \nabla_{X}Y\right)  ,Z\right)  $. On the other
hand, $\square$ has always a good restriction $\square:\mathfrak{X}\left(
\Sigma\right)  \times\mathfrak{X}\left(  \Sigma\right)  \rightarrow
\mathfrak{X}^{\ast}\left(  \Sigma\right)  $, which must coincide with
$\square^{\Sigma}$, the unique torsion-free metric dual connection on $\left(
\Sigma,g_{\Sigma}\right)  $. Since $\left(  \Sigma,g_{\Sigma}\right)  \,$is
riemannian, it must hold $\square_{X}^{\Sigma}Y\left(  Z\right)  =g_{\Sigma
}\left(  \nabla_{X}^{\Sigma}Y,Z\right)  $, and the result follows
$\clubsuit\vspace{0.5cm}$

The existence of a canonical radical vectorfield leads to the following Gauss formula.

\begin{proposition}
\label{2.0}Let $\left(  M,g\right)  $ be a transverse Riemann-Lorentz manifold
with transverse radical and $II$-flat. Then $\Sigma$ is \textquotedblright
totally geodesic\textquotedblright\ in the sense that, if $X,Y,Z,T\in
\mathfrak{X}\left(  \Sigma\right)  $ it holds
\[
K\left(  X,Y,Z,T\right)  =K^{\Sigma}\left(  X,Y,Z,T\right)
\]
where $K^{\Sigma}$ is the covariant curvature of $\Sigma$.
\end{proposition}

\textbf{Proof: }As we said in the proof of previous lemma we have, for
$X,Y,Z,T\in\mathfrak{X}\left(  \Sigma\right)  $: $\square_{X}Y\left(
Z\right)  =\square_{X}^{\Sigma}Y\left(  Z\right)  $, where $\square^{\Sigma} $
is the dual connection of $\left(  \Sigma,g_{\Sigma}\right)  $. Moreover,
since $\square_{X}R\left(  T\right)  =-\square_{X}T\left(  R\right)
=-II^{R}\left(  X,T\right)  =0$, again previous lemma leads to
\[
\square_{X}\left(  \nabla_{Y}Z\right)  \left(  T\right)  =\square_{X}\left(
\nabla_{Y}^{\Sigma}Z+III^{R}\left(  Y,Z\right)  R\right)  \left(  T\right)
=\square_{X}^{\Sigma}\left(  \nabla_{Y}^{\Sigma}Z\right)  \left(  T\right)
\]
what gives the result $\clubsuit\vspace{0.5cm}$

\begin{corollary}
Let $\left(  M,g\right)  $ be a transverse Riemann-Lorentz manifold with
transverse radical. If $\left(  \mathbb{M},g\right)  $ is flat, then $\left(
M,g\right)  $ is $III$-flat and $\Sigma$ is flat.
\end{corollary}

\textbf{Proof: }If $K=0$ then $Ric=0$. In particular, $Ric$ extends to $M$,
thus by Theorem \ref{1.1}, $\left(  M,g\right)  $ is $III$-flat. By
Proposition \ref{2.0}, $\Sigma$ is flat $\clubsuit\vspace{0.5cm}$

We now restate Theorem 9 of \cite{Kos87} in the following terms (the flatness
of $\Sigma$, being a consequence of the Collorary, needs not be included as an
extra hypothesis):

\begin{theorem}
\label{2.1}Let $\left(  M,g\right)  $ be a transverse Riemann-Lorentz
manifold. Then, $M$ is locally flat around $\Sigma$ if and only if, around
every singular point $p\in\Sigma$, there exists a coordinate system $\left(
\mathbb{U},x\right)  $ such that $g=\sum_{i=0}^{m-1}\left(  dx^{i}\right)
^{2}+\tau\left(  dx^{m}\right)  ^{2}$, where $\tau=0$ is a local equation for
$\Sigma$.
\end{theorem}

\section{Conformal geometry and the extendibility of Weyl curvature}

Let us consider a transverse Riemann-Lorentz manifold $\left(  M,g\right)  $
and the family $\mathcal{C}=\left\{  e^{2f}g:f\in C^{\infty}\left(  M\right)
\right\}  $. Take $\overline{g}=e^{2f}g\in\mathcal{C}$. Then $\left(
M,\overline{g}\right)  $ is also a transverse Riemann-Lorentz manifold, and
$\overline{\Sigma}=\Sigma$. Moreover, for each singular point $p\in\Sigma$ the
radical subspaces are the same: $\overline{Rad}_{p}\left(  M\right)
=Rad_{p}\left(  M\right)  $. We say that\emph{\ }$\left(  M,\mathcal{C}%
\right)  $\emph{\ is a transverse Riemann-Lorentz conformal manifold if some
(and thus any) }$g\in\mathcal{C}$\emph{\ is transverse Riemann-Lorentz}. Let
$\left(  M,\mathcal{C}\right)  $ be a transverse Riemann-Lorentz conformal
manifold. We say that $g\in\mathcal{C}$ \emph{is conformally }$II$\emph{-flat
if }$II_{\Sigma}^{R}=hg_{\Sigma}$\emph{, for some radical vectorfield }$R$
\emph{and some }$h\in C^{\infty}\left(  \Sigma\right)  $. This definition does
not depend on $R$ and, even more, it is conformal: if $\overline{g}=e^{2f}%
g\in\mathcal{C}$, then it holds
\begin{equation}
\overline{II}_{\Sigma}^{R}=e^{2f}\left\{  II_{\Sigma}^{R}-\left.  Rf\right|
_{\Sigma}g_{\Sigma}\right\}  \label{3.1}%
\end{equation}
Thus we say that $\left(  M,\mathcal{C}\right)  $ is \emph{conformally }%
$II$\emph{-flat if some (and thus, any) metric }$g\in\mathcal{C}$\emph{\ is
conformally }$II$\emph{-flat}.

\begin{proposition}
\label{3.2}A transverse Riemann-Lorentz conformal manifold $\left(
M,\mathcal{C}\right)  $ is conformally $II$-flat if and only if around every
singular point $p\in\Sigma$ there exist an open neighbourhood $\mathbb{U}$ in
$M$ and a metric $g\in\mathcal{C}$ which is $II$-flat on $\mathbb{U}$, that is
$II_{\Sigma\cap\mathbb{U}}=0$.
\end{proposition}

\textbf{Proof: }Let $\left(  \mathbb{U},E\right)  $ be an adapted orthonormal
frame near $p\in\Sigma$ (that is, $E_{m}$ is radical and $\left(
E_{1},...,E_{m-1}\right)  $ are orthonormal) and $g\in\mathcal{C}$. If
$\mathcal{C}$ is conformally $II$-flat, then there exists $h\in C^{\infty
}\left(  \Sigma\right)  $ such that $II_{\Sigma}^{E_{m}}=hg_{\Sigma}$. Take
$\widehat{h}\in C^{\infty}\left(  \mathbb{U}\right)  $ any local extension of
$h$ (shrinking $\mathbb{U}$ if necessary). There exists $f\in C^{\infty
}\left(  \mathbb{U}\right)  $ (shrinking again $\mathbb{U}$ if necessary)
satisfying $E_{m}f=\widehat{h}$ (since it is locally a first order linear
equation), what gives on $\mathbb{U}$: $II_{\Sigma}^{E_{m}}=\left.  \left(
E_{m}f\right)  \right\vert _{\Sigma}g_{\Sigma}$. Let $\widehat{f}\in
C^{\infty}\left(  M\right)  $ be any extension of (possibly a restriction of)
$f$. Applying (\ref{3.1}) to $g$ and $\overline{g}:=e^{2\widehat{f}}%
g\in\mathcal{C}$ we have $\overline{II}_{\Sigma}^{E_{m}}=0$.

To show the converse we start considering $g\in\mathcal{C}$. Since conformally
$II$-flatness is a local condition, it suffices to take an arbitrary
$p\in\Sigma$ and $\overline{g}=e^{2\widehat{f}}g\in\mathcal{C}$ such that
$\overline{g}$ is $II$-flat around $p$. Then, formula (\ref{3.1}) applied to
$g$ and $\overline{g}$ shows that $II_{p}^{\xi}=\left(  \xi f\right)  g_{p}$,
where $\xi\in Rad_{p}\left(  M\right)  -\left\{  0\right\}  $ $\clubsuit
\vspace{0.5cm}$

In what follows, we study \emph{conformally }$II$\emph{-flat Riemann-Lorentz
conformal} structures with transverse radical. Let $g$ and $\overline
{g}=e^{2f}g\in\mathcal{C}$ be two transverse Riemann-Lorentz metrics which are
$II$-flat. Formula (\ref{3.1}) shows that $\left.  \left(  Rf\right)  \right|
_{\Sigma}=0$. The expression of $grad_{g}\left(  f\right)  $ in an adapted
orthonormal frame such that $R=E_{m}$ is $grad_{g}\left(  f\right)
=\sum_{i=1}^{m-1}\left(  E_{i}f\right)  E_{i}+\tau^{-1}\left(  Rf\right)  R$,
thus $grad_{g}\left(  f\right)  $ extends to the whole $M$. Now a simple
computation gives
\begin{equation}
\overline{III}^{R}=e^{2f}\left\{  III^{R}-II^{R}\left(  grad_{g}\left(
f\right)  ,R\right)  g_{\Sigma}\right\}  \label{3.3}%
\end{equation}

We say that $g\in\mathcal{C}$ is \emph{conformally\ }$III$\emph{-flat}
\emph{if it is }$II$\emph{-flat (in order that }$III^{R}$\emph{\ exists) and
it holds }$III^{R}=kg_{\Sigma}$\emph{, for some radical vectorfield }%
$R$\emph{\ and some }$k\in C^{\infty}\left(  \Sigma\right)  $. Since
$II$-flatness is not conformal, the above definition, although independent of
$R$, cannot be conformal. However, it is conformal in the subset of $II$-flat
metrics$.$

\begin{definition}
We say that a transverse Riemann-Lorentz conformal manifold $\left(
M,\mathcal{C}\right)  $ with transverse radical is conformally $III$-flat if
it is conformally $II$-flat and every $g\in\mathcal{C}$ which is $II$-flat on
some open $\mathbb{U}$ of $M$ is also conformally $III$-flat on $\mathbb{U}$.
\end{definition}

Note that there may exist no conformally $III$-flat metrics on a conformally
$III$-flat manifold, simply because there may exist no $II$-flat metric there.
However, since a conformally $III$-flat space is conformally $II$-flat, we
deduce from Proposition \ref{3.2} that there always exist locally $II$-flat
metrics. Let us show that in fact there also exist locally $III$-flat metrics:

\begin{proposition}
\label{3.4}A transverse Riemann-Lorentz conformal manifold $\left(
M,\mathcal{C}\right)  $ with transverse radical is conformally $III$-flat if
and only if around every singular point $p\in\Sigma$ there exist an open
neighbourhood $\mathbb{U}$ in $M$ and a metric $g\in\mathcal{C}$ which is
$III$-flat on $\mathbb{U}$, that is $III_{\Sigma\cap\mathbb{U}}=0$.
\end{proposition}

\textbf{Proof: }Consider $p\in\Sigma$ and $\left(  \mathbb{U},E\right)  $ a
\emph{completely adapted orthonormal frame} (i.e., $E_{m}$ is radical and
$\left(  E_{1},...,E_{m-1}\right)  $ are orthonormal and tangent to $\Sigma$).
If $\left(  M,\mathcal{C}\right)  $ is conformally $III$-flat, there exist
$g\in\mathcal{C}$ which is $II$-flat on $\mathbb{U}$ (without loss of
generality) and $k\in C^{\infty}\left(  \Sigma\cap\mathbb{U}\right)  $, such
that $III^{E_{m}}=kg_{\Sigma}$. Since the radical is transverse, we have
$II_{mm}^{E_{m}}\neq0$, thus $k_{1}:=\frac{k}{II_{mm}^{E_{m}}}$ is $C^{\infty
}$ on $\Sigma\cap\mathbb{U}$. As in Proposition \ref{3.2} we can obtain $f\in
C^{\infty}\left(  \mathbb{U}\right)  $ such that $E_{m}f=\tau\widehat{k}_{1}$,
where $\tau=g\left(  E_{m},E_{m}\right)  $ and $\widehat{k}_{1}\in C^{\infty
}\left(  \mathbb{U}\right)  $ is any local extension of $k_{1}$. Since
$\left.  \left(  E_{m}f\right)  \right\vert _{\Sigma}=0$, we get
$grad_{g}\left(  f\right)  \in\mathfrak{X}\left(  \mathbb{U}\right)  $ and we
have $II^{E_{m}}\left(  grad_{g}\left(  f\right)  ,E_{m}\right)  =\left(
\tau^{-1}E_{m}f\right)  _{\Sigma}II_{mm}^{E_{m}}=k$. Now, take any extension
$\widehat{f}\in C^{\infty}\left(  M\right)  $ of (possibly a restriction of)
$f $. Since $g$ is $II$-flat, we deduce from (\ref{3.1}) that $\overline
{g}=e^{2\widehat{f}}g\in\mathcal{C}$ is also $II$-flat on $\mathbb{U}$. We
also deduce that $\overline{g}$ is $III$-flat on $\mathbb{U}$.

To prove the converse, first observe that the hypothesis implies in particular
that $\left(  M,\mathcal{C}\right)  $ is conformally $II$-flat. Consider
$p\in\Sigma$ and $g\in\mathcal{C}$, $II$-flat on a neighbourhood of $p$. By
hypothesis, there exists $\overline{g}=e^{2f}g\in\mathcal{C}$ which is
$III$-flat around $p$. Thus we deduce from (\ref{3.3}) that $III^{R}%
=II^{R}\left(  grad_{g}\left(  f\right)  ,R\right)  g_{\Sigma}$, so $g$ is
conformally $III$-flat $\clubsuit\vspace{0.5cm}$

In what follows we shall assume that $\dim M=m\geq4$. We now study the
extendibility of $\emph{the}$ \emph{Weyl tensor}, naturally defined on
$\left(  \mathbb{M},\mathcal{C}_{\mathbb{M}}\right)  $. It is well-known that
this tensor plays a main role in deciding when $\mathbb{M}$ is (locally)
conformally flat, according to \emph{Weyl Theorem: a pseudoriemannian
conformal manifold is (locally) conformally flat if and only if the Weyl
tensor vanishes identically }(see for instance the preliminaries of
\cite{He}). At the end of the paper we discuss the problem of establish a
modified version of Weyl Theorem for transverse Riemann-Lorentz conformal manifolds.

The \emph{Weyl tensor }$W$ on $\left(  \mathbb{M},g_{\mathbb{M}}\right)  $ can
be defined as
\[
W:=K-h\bullet g\in\mathcal{I}_{4}^{0}\left(  \mathbb{M}\right)  ,
\]
where $h=\frac{1}{m-2}\left\{  Ric-\frac{Sc}{2\left(  m-1\right)  }g\right\}
$ is the \emph{Schouten tensor}, $Ric$ is the Ricci tensor and $Sc$ is the
scalar curvature associated to $\left(  \mathbb{M},g_{\mathbb{M}}\right)  $,
and where
\[
\bullet:S^{2}\left(  \mathbb{M}\right)  \times S^{2}\left(  \mathbb{M}\right)
\rightarrow\mathcal{I}_{4}^{0}\left(  \mathbb{M}\right)
\]
is the so-called \emph{Kulkarni-Nomizu product}, given by
\[
\theta\bullet\omega\left(  x,y,z,t\right)  :=\det\left(
\begin{array}
[c]{ll}%
\theta\left(  x,z\right)  & \omega\left(  x,t\right) \\
\theta\left(  y,z\right)  & \omega\left(  y,t\right)
\end{array}
\right)  +\det\left(
\begin{array}
[c]{ll}%
\omega\left(  x,z\right)  & \theta\left(  x,t\right) \\
\omega\left(  y,z\right)  & \theta\left(  y,t\right)
\end{array}
\right)
\]
If we pick $\overline{g}=e^{2f}g\in\mathcal{C}$, then the Weyl tensor
associated to $\left(  \mathbb{M},\overline{g}_{\mathbb{M}}\right)  $
satisfies $\overline{W}=e^{2f}W$, thus \emph{the Weyl conformal curvature
}$\mathcal{W}:=\uparrow_{2}^{1}W\in\mathcal{I}_{3}^{1}\left(  \mathbb{M}%
\right)  $ becomes a conformal invariant. Notice that the extendibility of $W$
(which is equivalent to the extendibility of $\mathcal{W}$) is a conformal
condition, therefore it should be stated in terms of the conformal structure.
In fact, we prove that it is equivalent to conformal $III$-flatness.

\begin{theorem}
\label{3.5}Let $\left(  M,\mathcal{C}\right)  $ be a transverse
Riemann-Lorentz conformal manifold, with $\dim M=m\geq4$. Then $W$ (smoothly)
extends to the whole $M$ if and only if the radical is transverse and
$\mathcal{C}$ is conformally $III$-flat.
\end{theorem}

\textbf{Proof: }If $\left(  M,\mathcal{C}\right)  $ has transverse radical and
is conformally $III$-flat, there exist (Proposition \ref{3.4}) a $M$-open
covering $\left\{  \mathbb{U}_{\alpha}\right\}  $ of $\Sigma$ and a family of
metrics $\left\{  g_{\alpha}\right\}  $ in $\mathcal{C}$ such that $g_{\alpha
}$ is $III$-flat on $\mathbb{U}_{\alpha}$. By Theorem \ref{1.1}, the covariant
curvature $K_{\alpha}$, the Ricci tensor $Ric_{\alpha}$ and the scalar
curvature $Sc_{\alpha}$ associated to $g_{\alpha}$ extend to $\Sigma
\cap\mathbb{U}_{\alpha}$, therefore the Weyl tensor $W_{\alpha}$ also extends
to $\Sigma\cap\mathbb{U}_{\alpha}$. Since this is a conformal condition,
$W_{\alpha}$ extends to $\Sigma\cap\mathbb{U}_{\beta}$ for all $\beta$, and
thus $W_{\alpha}$ extends to the whole $M$.

To show the converse we start picking an adapted orthonormal frame $\left(
\mathbb{U},E\right)  $. Then, we can express the functions $W_{abcd}=W\left(
E_{a},E_{b},E_{c},E_{d}\right)  $ as second order polynomials in $\tau
^{-1}=\left(  g\left(  E_{m},E_{m}\right)  \right)  $. Let us call $\left(
W_{abcd}\right)  _{0},\left(  W_{abcd}\right)  _{1}$ and $\left(
W_{abcd}\right)  _{2}$ the differentiable coefficients of the terms of order
$0$, $1$ and $2$. Since $\tau=0$ is a local equation for $\Sigma$, $W$ extends
to $\mathbb{U}$ if and only if the restricted functions $\left.  \left(
W_{abcd}\right)  _{2}\right|  _{\Sigma}$ and $\left.  \left(  W_{abcd}\right)
_{1}+\tau^{-1}\left(  W_{abcd}\right)  _{2}\right|  _{\Sigma}$ identically vanish.

Suposse the radical is tangent to $\Sigma$ at a singular point $p\in\Sigma$.
We can choose the frame such that $E_{1}\left(  p\right)  ,E_{2}\left(
p\right)  \in T_{p}M-T_{p}\Sigma$. But then, using that $II^{E_{m}}\left(
E_{m},E_{m}\right)  \left(  p\right)  =0$ (because the radical is tangent), we
obtain $\left(  W_{1323}\left(  p\right)  \right)  _{2}=\frac{\varepsilon_{3}%
}{m-2}II_{p}^{E_{m}}\left(  E_{1},E_{m}\right)  II_{p}^{E_{m}}\left(
E_{2},E_{m}\right)  $. Since $E_{1}$ and $E_{2}$ are transverse to $\Sigma$ at
$p$, $\left(  W_{1323}\left(  p\right)  \right)  _{2}\neq0$, hence $W$ cannot
be extended. Therefore the radical must be transverse to $\Sigma$.

Once we know that the radical must be always transverse to $\Sigma$ (thus
$II_{mm}^{E_{m}}\neq0$), we can choose the orthonormal frame $\left(
\mathbb{U},E\right)  $ completely adapted. Thus, picking $i,j,k$ different
from $m$, with $i,j$ different from $k$, and using that $II_{im}^{E_{m}}=0$,
we have: if $i\neq j$, then $0=\left.  \left(  W_{ikjk}\right)  _{2}\right|
_{\Sigma}=-\frac{\varepsilon_{k}}{m-2}II_{ij}^{E_{m}}II_{mm}^{E_{m}}$. Since
$II_{mm}^{E_{m}}\neq0$, we get $II_{ij}^{E_{m}}=0$. If $i=j$ (and using
$II_{ij}^{E_{m}}=0$), the $\left(
\begin{array}
[c]{c}%
m-1\\
2
\end{array}
\right)  $ equalities $0=\left.  \left(  W_{ikik}\right)  _{2}\right|
_{\Sigma}$, suitably manipulated, give us $\varepsilon_{i}II_{ii}^{E_{m}%
}+\varepsilon_{k}II_{kk}^{E_{m}}=\frac{2C}{m-1}$, where $C=\sum_{l=1}%
^{m-1}\varepsilon_{l}II_{ll}^{E_{m}}\in C^{\infty}\left(  \mathbb{U}\right)
$. Substracting the equation for $i,k$ from the equation for $k,j$, we obtain
$\varepsilon_{i}II_{ii}^{E_{m}}-\varepsilon_{j}II_{jj}^{E_{m}}=0$, thus
$\varepsilon_{i}II_{ii}^{E_{m}}=\varepsilon_{j}II_{jj}^{E_{m}}$. Defining
$k:=\varepsilon_{1}II_{11}^{E_{m}}\in C^{\infty}\left(  \Sigma\cap
\mathbb{U}\right)  $, it holds $II_{ii}^{E_{m}}=\varepsilon_{i}\varepsilon
_{1}II_{11}^{E_{m}}=kg_{ii}$ and $II_{ij}^{E_{m}}=0=kg_{ij}$ (where $i\neq
j$), what means $II_{\Sigma}^{E_{m}}=kg_{\Sigma}$, that is, $g$ is conformally
$II$-flat on $\mathbb{U}$, and therefore $\left(  M,\mathcal{C}\right)  $ is
conformally $II$-flat.

Once we know that $\left(  M,\mathcal{C}\right)  $ is conformally $II$-flat,
we can choose a metric $g\in\mathcal{C}$ which is $II$-flat on $\mathbb{U}$
(shrinking $\mathbb{U}$ if neccesary). By Theorem \ref{1.1}, the covariant
curvature $K$ associated to $g$ extends to $\Sigma\cap\mathbb{U}$ and, since
$W$ also does it, necessarily $h\bullet g$ extends to $\Sigma\cap\mathbb{U}$.
Picking $i,j,k$ different from $m$, with $i,j$ different from $k $, we get
$\left(  h\bullet g\right)  _{ikjk}=\varepsilon_{k}h_{ij}+\delta
_{ij}\varepsilon_{i}h_{kk}=A_{ijk}+\tau^{-1}B_{ijk}$, therefore the function
\[
B_{ijk}:=\frac{1}{m-2}\left\{  \varepsilon_{k}K_{imjm}+\delta_{ij}%
\varepsilon_{i}K_{kmkm}-\frac{2\varepsilon_{k}\delta_{ij}\varepsilon_{i}}%
{m-1}\sum_{l=1}^{m-1}\varepsilon_{l}K_{lmlm}\right\}
\]
must vanish on $\Sigma$. Using the same argument as before, but with the
equalities $0=\left.  B_{ijk}\right\vert _{\Sigma}$, we get $III^{E_{m}%
}=kg_{\Sigma}$, where $k:=\varepsilon_{1}III_{11}^{E_{m}}\in C^{\infty}\left(
\Sigma\cap\mathbb{U}\right)  $, that is $g$ is conformally $III$-flat on
$\mathbb{U}$, and thus $\left(  M,\mathcal{C}\right)  $ is conformally
$III$-flat $\clubsuit\vspace{0.5cm}$

Let us consider the following conjecture:

\begin{conjecture}
\label{conjecture}Let $\left(  M,\mathcal{C}\right)  $ be a transverse
Riemann-Lorentz conformal manifold, with $\dim M=m\geq4$. A necessary
condition for being $W=0$ is that, around every singular point $p\in\Sigma$,
there exist a coordinate system $\left(  \mathbb{U},x\right)  $ and a metric
$g\in\mathcal{C}$ such that $g=\sum_{i=0}^{m-1}\left(  dx^{i}\right)
^{2}+\tau\left(  dx^{m}\right)  ^{2}$, where $\tau=0$ is a local equation for
$\Sigma$.
\end{conjecture}

Using Theorem \ref{2.1}, it becomes obvious that the necessary condition
stated in the conjecture is always sufficient for having $W=0$ around $\Sigma$.

If the conjecture is true, $\Sigma$ must be (locally) conformally flat, which
is well known equivalent to either $W^{\Sigma}=0$ (if $m>4$) or $\nabla
_{X}^{\Sigma}h^{\Sigma}(Y,Z)=\nabla_{Y}^{\Sigma}h^{\Sigma}(X,Z)$ (if $m=4$).
But the extendibility of $W$, equivalent (Theorem \ref{3.5}) to conformal
$III$-flatness, implies (Proposition \ref{3.4}) the existence of a metric
$g\in\mathcal{C}$ which is $III$-flat around $\Sigma$, thus satisfying
(Proposition \ref{2.0}):%
\[
W\mid_{T\Sigma}\quad=\left(  K-h\bullet g\right)  \mid_{T\Sigma}%
\quad=K^{\Sigma}-h\mid_{T\Sigma}\bullet g_{\Sigma}=W^{\Sigma}+(h^{\Sigma
}-h\mid_{T\Sigma})\bullet g_{\Sigma}\quad.
\]
Because conditions $W=0$ and $W^{\Sigma}=0$ are conformal, any counterexample
$\left(  M,\mathcal{C}\right)  $ to the above conjecture must admit a metric
$g\in\mathcal{C}$ which is $III$-flat around $\Sigma$ and satisfies either
$h^{\Sigma}\neq h\mid_{T\Sigma}$ (if $m>4$) or (Lemma \ref{2.0.0}) $\nabla
_{X}h(Y,Z)\neq\nabla_{Y}h(X,Z)$, for some $X,Y,Z\in\mathfrak{X}(\Sigma)$ (if
$m=4$). Now a straightforward computation for $III $-flat metrics, using an
orthonormal completely adapted frame, leads to the following expression in
terms of extendible quantities:%
\[
h_{ij}^{\Sigma}-h_{ij}\mid_{T\Sigma}\quad=\frac{-1}{m-2}\left\{
\frac{K_{imjm}}{\tau}-\frac{1}{m-3}\sum_{l=1}^{m-1}K_{iljl}-\right.
\quad\quad\quad\quad
\]%
\[
\quad\quad\quad\quad\quad\quad\quad\quad\quad\quad\left.  -\frac{1}%
{m-1}\left[  \sum_{k=1}^{m-1}\frac{K_{kmkm}}{\tau}-\frac{1}{m-3}\sum
_{k,l=1}^{m-1}K_{klkl}\right]  \delta_{ij}\right\}  \mid_{\Sigma}\quad,
\]
($i,j=1,...,m-1$), which shows that the construction of counterexamples is not
easy.\vspace{0.5cm}

In fact, \emph{the conjecture is true for transverse Riemann-Lorentz warped
products}, as we show right now. Let us consider a $m$-dimensional ($m\geq4$)
transverse Riemann-Lorentz manifold $\left(  M,g\right)  $ of the form
$M=I\times S$, where $\dim I=1$, $0\in I$, and $g=f\left(  t\right)  ^{2}%
g_{S}-tdt^{2}$, where $f\in C^{\infty}\left(  I\right)  $, $f>0$ and $g_{S}$
is riemannian (we identify $t$, $f$ and $g_{S}$ with the corresponding
pullbacks by the canonical projections). Thus $\Sigma=\left\{  0\right\}
\times S$ is homothetic to $S$ with scale factor $f(0)$. Calling
$U\in\mathfrak{X}\left(  M\right)  $ the (nowhere zero) lift of the
vectorfield $\frac{d}{dt}\in\mathfrak{X}\left(  I\right)  $, one inmediately
sees that $U$ is radical and transverse to $\Sigma$. It is not difficult to
compute the curvature tensors on $\mathbb{M}$. Standar results on warped
products (see \cite{O'Ne}, Chapter 7) lead to (we denote by $X,Y\in
\mathfrak{X}\left(  M\right)  $ the lifts of corresponding vectorfields
$\overline{X},\overline{Y}\in\mathfrak{X}\left(  S\right)  $) $\nabla
_{U}U=\frac{1}{2t}U$, $\nabla_{U}X=\nabla_{X}U=\frac{f^{\prime}}{f}X$ and
$\nabla_{X}Y=g\left(  X,Y\right)  \frac{f^{\prime}}{tf}U+\nabla_{\overline{X}%
}^{S}\overline{Y}$ (where $\nabla^{S}$ is the Levi-Civita connection on $S$
and $\nabla_{\overline{X}}^{S}\overline{Y}$ is the lift of the corresponding
vectorfield on $S$) and also the following expressions for the curvature
tensors:%
\[
\left\{
\begin{array}
[c]{l}%
K=f^{2}K^{S}+\frac{f^{\prime2}f^{2}}{2t}g_{S}\bullet g_{S}+\frac{f}{2}%
(\frac{f^{\prime}}{t}-2f^{\prime\prime})g_{S}\bullet dt^{2}\\
Ric=Ric^{S}-\left(  \frac{f}{2t}(\frac{f^{\prime}}{t}-2f^{\prime\prime
})-\left(  m-2\right)  \frac{f^{\prime2}}{t}\right)  g_{S}+\frac{m-1}%
{2f}(\frac{f^{\prime}}{t}-2f^{\prime\prime})dt^{2}\\
Sc=\frac{Sc^{S}}{f^{2}}-\frac{m-1}{f^{2}}\left(  \frac{f}{t}(\frac{f^{\prime}%
}{t}-2f^{\prime\prime})-\left(  m-2\right)  \frac{f^{\prime2}}{t}\right) \\
h=\frac{m-3}{m-2}h^{S}+\left(  \frac{Sc^{S}}{2(m-2)^{2}\left(  m-1\right)
}+\frac{f^{\prime2}}{2t}\right)  g_{S}+\\
\quad\quad\quad\quad\quad\quad\quad\quad\quad+\left(  \frac{tSc^{S}}{2\left(
m-1\right)  \left(  m-2\right)  f^{2}}+\frac{1}{2f}(\frac{f^{\prime2}}%
{f}+\frac{f^{\prime}}{t}-2f^{\prime\prime})\right)  dt^{2}\\
W=f^{2}W^{S}+\frac{1}{(m-2)}\left(  Ric^{S}-\frac{Sc^{S}}{m-1}g_{S}\right)
\bullet\left(  \frac{f^{2}}{m-3}g_{S}+tdt^{2}\right)
\end{array}
\right.
\]
($K^{S}$, $Ric^{S}$, $Sc^{S}$, $h^{S}$ and $W^{S}$ denote of course the
pullbacks by the projection of the corresponding tensor fields on $S$). It follows:

\begin{lemma}
The following three conditions are equivalent: (1) $K$ extends to $M$, (2)
$f^{\prime}\left(  0\right)  =0$ and (3) $h$ extends to $M$. Also the
following are equivalent: (1) $Ric$ extends to $M$, (2) $\left(  f^{\prime
}/t\right)  \left(  0\right)  =0$ and (3) $Sc$ extends to $M$. Moreover, $W$
extends to $M$ in any case.
\end{lemma}

The fact that $W$ extends to $M$ was obvious from the very beginning: the map
$\Psi\equiv\psi\times id:\left(  I-\left\{  0\right\}  \right)  \times
S\rightarrow\mathbb{R}\times S$, given by $T\equiv\psi\left(  t\right)
:=\int_{0}^{t}\frac{\left\vert s\right\vert ^{\frac{1}{2}}ds}{f\left(
s\right)  }$, is a conformal diffeomorphism onto its (non-connected) image
with the metric $\overline{g}\equiv-\left(  dT\right)  ^{2}+g_{S}$, thus it
preserves the $(_{3}^{1})$-Weyl tensors, and since $\overline{g}$ is regular
around $T=0$ and $f\left(  0\right)  \neq0$, $\overline{W}$ (and therefore
$W$) extends to the whole $M$. It follows from Theorem \ref{3.5} that the
conformal manifold $\left(  M,\left[  g\right]  \right)  $ is (in any case)
conformally $III$-flat.

\begin{lemma}
\label{3.6}The following four conditions are equivalent: (1) $W=0$, (2)
$W^{S}=0=Ric^{S}-\frac{Sc^{S}}{m-1}g_{S}$ and (3) $\Sigma$ has constant
(sectional) curvature.
\end{lemma}

\textbf{Proof: }$(1)\Leftrightarrow(2)$ follows from the above formula.
$(2)\Rightarrow(3)$: $Ric^{S}-\frac{Sc^{S}}{m-1}g_{S}=0$ implies (Schur's
lemma) $Sc^{S}=(m-1)(m-2)C$ (constant), thus $h^{S}=\frac{C}{2}g_{S}$;
moreover $W^{S}=0$ leads to $K^{S}=\frac{C}{2}g_{S}\bullet g_{S}$.
$(3)\Rightarrow(2)$: From $K^{S}=\frac{C}{2}g_{S}\bullet g_{S}$, one
immediately gets $W^{S}=0=Ric^{S}-\frac{Sc^{S}}{m-1}g_{S}$ $\clubsuit$

\begin{proposition}
The Conjecture \ref{conjecture} is true for any transverse Riemann-Lorentz
conformal manifold $\left(  M,\mathcal{C}\right)  $ such that some
$g\in\mathcal{C}$ is a warped product.
\end{proposition}

\textbf{Proof: }Let $g=f\left(  t\right)  ^{2}g_{S}-tdt^{2}\in\mathcal{C}$ be
a transverse warped product metric on $M=I\times S$. Note that $g=f\left(
t\right)  ^{2}\left\{  g_{S}-\frac{t}{f\left(  t\right)  ^{2}}dt^{2}\right\}
$. From $W=0$ and Lemma \ref{3.6} we get, around any $p\in\Sigma$, coordinates
$\left(  \mathbb{V},y\right)  $ of $\Sigma$ such that $f\left(  0\right)
^{2}g_{S}=g_{\Sigma}=e^{2h}\sum_{i=1}^{m-1}\left(  dy^{i}\right)  ^{2}$, for
some $h\in C^{\infty}(\Sigma)$. Choosing $x^{i}:=y^{i}\circ\pi$, $x^{m}:=t$
and $\tau:=\frac{-te^{-2h}}{f\left(  t\right)  ^{2}}$, we get $g=e^{2h}%
f\left(  t\right)  ^{2}\left\{  \sum_{i=1}^{m-1}\left(  dx^{i}\right)
^{2}+\tau\left(  dx^{m}\right)  ^{2}\right\}  $, and we are finished
$\clubsuit$


\begin{thebibliography}{9}                                                                                                %


\bibitem {Ag-La}E. Aguirre and J. Lafuente. \emph{Trasverse Riemann-Lorentz
metrics with tangent radical}. Diff. Geom. its App, 24, 2, 91-100, 2005.

\bibitem {Ha-Ha}J. B. Hartle and S. W. Hawking. \emph{Wave Function of the
Universe}. Phys. Rev., D41, 1815-34, 1990.

\bibitem {He}U. Hertrich-Jeromin. \emph{Introduction to M\"{o}bius
Differential Geometry.} Cambridge Univ. Press, 2003.

\bibitem {Kos85}M. Kossowski. \emph{Fold singularities in pseudoriemannian
geodesic tubes. }Proc. Amer. Math. Soc., 95, 463-469,1985.

\bibitem {Kos87}M. Kossowski. \emph{Pseudo-riemannian metric singularities and
the extendability of parallel transport}. Proc. Amer. Math. Soc., 99, 147-154, 1987.

\bibitem {Kos94}M. Kossowski and M. Kriele. \emph{Transverse, type changing,
pseudo riemannian metrics and the extendability of geodesics. }Proc. R. Soc.
Lond. A 444, 297-306,1994.

\bibitem {Kos97}M. Kossowski and M. Kriele. \emph{The volume blow-up and
characteristic classes for transverse, type changing, pseudo-riemannian
metrics.} Geom. Dedicata 64, 1-16. 1997.

\bibitem {O'Ne}B. O`Neill. \emph{Semi-riemannian Geometry.} Academic Press, 1983.
\end{thebibliography}
\end{document}